%% file: AdaptIV.tex
\documentclass[12pt,a4paper,final]{article}

\usepackage[a4paper,portrait]{geometry}
\usepackage{amsmath,amsfonts,latexsym,amssymb,amsthm}
\usepackage{graphicx}

\newcommand{\mykeywords}{
 Inverse Problems,
Instrumental Variables,
Model Selection,
Econometrics
.}

\newcommand{\mysubjclass}{
62G05, 62G20
.}
\input{macros}

\THMEN
\title{Oracle Inequality for Instrumental Variable Regression}
\author{J-M. Loubes \& C. Marteau}
\begin{document}
\maketitle \noindent
\begin{abstract}
We tackle the problem of estimating a regression function observed in an instrumental regression framework. This model is an inverse problem with unknown operator. We provide a spectral cut-off estimation procedure which enables to derive oracle inequalities which warrants that our estimate, built without any prior knowledge, behaves as well as, up to $\log$ term, if the best model were known.
\end{abstract}

{ \noindent
 \textbf{Subject Class. MSC-2000}: \mysubjclass}

\section*{Introduction}
An economic relationship between a response variable $Y$ and a vector of explanatory variables $X$ is often represented by an equation
$$Y=\vphi(X)+U,$$ where $\vphi$ is the parameter of interest which models the relationship while $U$ is an error term. Contrary to usual statistical regression models, the error term is correlated with the explanatory variables $X$, hence $\bE (U|X)\neq 0$, preventing direct estimation of $\vphi$. To overcome the endogeneity of $X$, we assume that there exists   an observed random variable $W$, called the instrument, which decorrelates the effects of the two variables $X$ and $Y$ in the sense that $\bE(U|W)=0$. It is often the case in economics, where the practical construction of instrumental variables play an important part.  For instance \cite{MR2380756} present practical situations where prices of goods and quantity in goods can be explained using an instrument. This situation is also encountered when dealing with simultaneous equations, error-in-variable models, treatment model with endogenous effects. It defines the so-called instrumental variable regression model which has received a growing interest among the last decade and turned to be a challenging issue in statistics.  In particular, we refer to \cite{hausnew}, \cite{MR2000257} \cite{handbook} for general references on the use of instrumental variables in economics while \cite{MR2253107}, \cite{florens} and \cite{joflo} deal with the statistical estimation problem.\vskip .1in
\vskip .1in
 More precisely, we aim at estimating a function $\vphi$ observed in the following observation model
\begin{equation} \label{model:IV}
 Y  =\vphi(X)+U,\: \quad \begin{cases} \bE(U|X) & \neq 0 \\
\bE(U|W) & = 0\end{cases}
\end{equation}
Hence, the model \eqref{model:IV} can be rewritten as an inverse problem using  the expectation conditional operator with respect to $W$, which will be denoted $T$,  as follows :
 \begin{equation} \label{model:INV}
 r :=\bE(Y|W)= \bE(\vphi(X)|W)=T\vphi .\end{equation}
The function $r$ is not known and only an observation $\hat{r}$ is available, leading to the inverse problem  $\hat{r}=T\vphi+\delta$, where $\vphi$ is defined as the solution of a noisy Fredholm equation of the first order which may generate an ill-posed inverse problem. The literature on inverse problems in statistics is large, but  contrary to most of the problems tackled in the literature on inverse problems (see \cite{engl96}, \cite{Ruym96}, \cite{cgpt}, \cite{galerkin}, \cite{loublud1} and \cite{osulliv} for general references), the operator $T$ is unknown either, which transforms the model into an inverse problem with unknown operator.   Few results exist in this settings and only very recently  new methods have arised. In particular   \cite{MR2158113}, \cite{clem,clem2}, or \cite{MR1872847} and \cite{MR2387973} in a more general case, construct estimators which enable to estimate inverse problem with {\it unobserved} operators in an adaptive way, i.e getting optimal rates of convergence without prior knowledge of the regularity of the functional parameter of interest.\vskip .1in

In this work, we are facing an even more difficult situation since both $r$ and the operator $T$ have to be estimated from the same sample.  Some attention has been paid to this estimation issue, with different kinds of technics such as kernel based Tikhonov regularization \cite{florens} or \cite{MR2253107}, regularization in Hilbert scales \cite{joflo}, finite dimensional sieve minimum distance estimator \cite{MR2000257}, with different rates and different smoothness assumptions,  providing sometimes minimax rates of convergence. But, to our knowledge, all the proposed estimators rely on  prior knowledge on the regularity of the function $\vphi$ expressed through an embedding condition into a smoothness space or an Hilbert scale, or a condition linking the regularity of $\vphi$ to the regularity of the operator, namely a link condition or source condition (see \cite{xiaoreiss} for general comments and insightful comments on such assumptions).\vskip .1in  Hence, in this paper, we provide under some conditions,  an adaptive estimation procedure of  the function $\vphi$ which converges, without prior regularity assumption, at the optimal rate of convergence, up to a logarithmic term. Moreover, we derive an oracle inequality which ensures optimality among the different choices of estimators.\vskip .1in
The article falls into the following parts. Section \ref{s:model} is devoted to the mathematical presentation of the instrumental variable framework and the building of the estimator. Section \ref{s:main} provides the asymptotic behaviour of this adaptive estimate as well as an oracle inequality, while technical Lemmas and proofs are gathered in Section \ref{s:appendix}.  
\section{Inverse Problem for IV regression} \label{s:model}
We observe  an i.i.d sample $(Y_i,X_i,W_i)$ for $i=1,\dots,n$ with unknown distribution $f(Y,X,W)$. Define the following Hilbert spaces 
\begin{align*}
L^2_X & =\{h: \mathbb{R}^d\to \bR, \: \|h\|^2_X:=\bE(h^2(X)) < + \infty \} \\
L^2_W & =\{g: \mathbb{R}^d\to \bR, \: \|g\|^2_W:=\bE(g^2(W)) < + \infty \}, 
\end{align*}
with the corresponding scalar product $<.,.>_X$ and $<.,.>_W$. Then the conditional expectation operator of $X$ with respect to $W$ is defined as an operator $T$
\begin{align*}
 T:\quad  L^2_X & \rightarrow L^2_W \\
 g & \rightarrow \bE (g(X)|W) . \end{align*}   
 The model \eqref{model:IV} can be written, as discussed in \cite{xiaoreiss}, as
 \begin{align} 
Y_i & = \vphi(X_i)+ \bE [\vphi(X_i)|W_i] - \bE [\vphi(X_i)|W_i] +U_i \nonumber   \\ 
& = \bE [\vphi(X_i)|W_i] +V_i  \nonumber \\ 
& = T\vphi(W_i)+V_i, \label{eq:modelreiss}
\end{align}
where $V_i= \vphi(X_i)-\bE [\vphi(X_i)|W_i] +U_i$, is such that $\bE (V|W)=0$. The parameter of interest is the unknown function $\vphi$. Hence, the observation model turns to be an inverse problem with unknown operator $T$ with a correlated noise $V$. Solving this issue amounts to deal with the estimation of the operator and then controlling the correlation with respect to the noise.
\vskip .1in
The operator $T$ is unknown and depends on the unknown distribution of the observed variables $f_{(Y,X,Z)}$. Estimation of an operator can be performed either by directly using an estimate of $f_{(Y,X,Z)}$, or if exists, by estimating the spectral value decomposition of the operator. \\
\indent Assume that $T$ is compact and admits a singular value decomposition (SVD) $(\la_j,\phi_j,\psi_j)_{j\geq 1}$, which provides a natural basis adapted to the operator for representing the function $\vphi$, see for instance \cite{engl96}. More precisely,  let $T^*$ be the adjoint operator of $T$, then $T^*T$ is a compact operator on $L^2_X$ with eigenvalues $\la_j^2,\: j\geq 1$ associated to the corresponding eigenfunctions $\phi_j$, while $\psi_j$ are defined by $\psi_j=\frac{T \phi_j }{\|T\phi_j\|}$. So we obtain
$$ T\phi_j=\la_j \psi_j,\quad T^* \psi_j=\la_j \phi_j.$$  We can write the following decompositions
\begin{equation} \label{est:eigen}  
r(w)=  \bE(Y|W=w)=T \vphi(w)   = \sum_{j\geq 1} \la_j <\vphi,\phi_j> \psi_j(w),  \end{equation}
\begin{equation}{\rm and} \quad r(w) = \sum_{j \geq 1} r_j \psi_j(w), \label{est:eigen2}
\end{equation}
with $r_j=<Y,\psi_j>$ that can be estimated by $$\hat{r}_j=\frac{1}{n}\sum_{i=1}^n Y_i \psi_j(W_i). $$ Hence the noisy observations are the $\hat{r}_j$'s which will be used to estimate the regression function $\vphi$ in an inverse problem framework.\vskip .1in 
In a very general framework, full estimation of an operator is a hard task  hence we restrict ourselves to the case where the SVD of the operator is partially known in the sense that the eigenvalues $\la_j$'s are unknown but the eigenvectors $\phi_j$'s and $\psi_j$'s are known.\\
\indent Note that this assumption is often met for the special case of deconvolution.
Consider the case where the unknown function $\vphi$ reduces to the identity. Hence model \eqref{model:IV} reduces to the usual deconvolution model
$$ Y=X +U .$$
Set $f_U$ the unknown density of the noise $U$ and assume that $f_U \in L^2(\bR)$ is a 1 periodic function. Let also $T_U$ be the convolution operator defined by $ T_U g = g \star f_U.$ In this special case, the spectral decomposition of the operator $T_U$ is known, given by the unitary Fourier transform and the usual real trigonometric basis on $[0,1]$ are the eigenvectors .\vskip .1in
If the operator were known we could provide an estimator using the spectral decomposition of the function $\vphi$ as follows. For a given decomposition level $m$, define the projection estimator (also called spectral cut-off \cite{engl96}) 
\begin{equation} \label{naif} \hat{\vphi}_m^0=\sum_{j=1}^m \frac{\hat{r}_j}{\la_j} \phi_j \end{equation}
Since the $\la_j$'s are unknown, we first build an estimator of the eigenvalues. For this, using the decomposition \eqref{est:eigen}, we obtain 
\begin{align*}
\la_j & =<T \phi_j,\psi_j>_W \\
 & =  \bE [T\phi_j(W) \psi_j(W)] \\
& = \bE [ \bE[\phi_j(X)|W] \psi_j(W)]\\
& = \bE [\phi_j(X) \psi_j(W)].
\end{align*}
So the eigenvalue $\la_j$ can be estimated by \begin{equation} \label{est:eig}
\hat{\lambda}_j=\frac{1}{n} \sum_{j=1}^n \psi_j(W_i) \phi_j(X_i). \end{equation}
As studied in \cite{MR2158113}, replacing directly the eigenvalues by their estimates in \eqref{naif} does not yield a consistent estimator, hence using their same strategy we define an upper bound for the resolution level
\begin{equation}
M= \inf \left\lbrace k\leq N: |\hat\lambda_k | \leq \frac{1}{\sqrt{n}} \log n \right\rbrace -1,
\label{M}
\end{equation}
for $N$  to be chosen later.  The parameter $N$ provides an upper bound for $M$ in order to ensure that $M$ is not too large.  The main idea behind this definition is that when the estimates of the eigenvalues are too small with respect to the observation noise, trying to still provide an estimation of the inverse $\la_k^{-1}$ only amplificates the estimation error. To avoid this trouble, we truncate the sequence of the estimated eigenvalues when their estimate is too small, i.e smaller than the noise level.  We point out that this parameter $M$ is a random variable which we will have to control. More precisely, define two deterministic lower and upper bounds
 $M_0,M_1$ as
\begin{equation}
M_0= \inf \left\lbrace k: |\lambda_k | \leq \frac{1}{\sqrt{n}} \log^2 n \right\rbrace -1,
\label{defM_0}
\end{equation}
and
\begin{equation}
M_1= \inf \left\lbrace k: |\lambda_k | \leq \frac{1}{\sqrt{n}} \log^{3/4} n \right\rbrace,
\label{defM_1}
\end{equation}
we will show in Section~\ref{s:appendix}, that with high probability $M_0 \leq M < M_1$. \vskip .1in
Now, thresholding the spectral decomposition in \eqref{naif} leads to the following estimator
\begin{equation}
	\hat{\vphi}_m=\sum_{j=1}^m  \frac{\hat{r}_j}{\hat{\la}_j}1_{j\leq M}\phi_j.
	\label{eq:estimateur}
\end{equation}
The asymptotic behaviour of this estimate depends on the choice of $m$. In the next section, we provide an optimal procedure to select the parameter $m$ that gives rise to an adaptive estimator $\vphi^\star$ and an oracle inequality.

\section{Main result}\label{s:main}
Consider the following assumptions on both the data $Y_i,\: i=1,\dots,n$ and the eigenfunctions $\phi_k$ and $\psi_k$ for $k\geq 1$. 
\begin{description}
	\item[Bounded SVD functions:] 
There exists a finite constant $C_1$ such that
\begin{equation} \label{basis}
\forall j\geq 1,\quad \|\phi_j\|_\infty < C_1, \quad  \|\psi_j\|_\infty< C_1 \end{equation}
\end{description}
\begin{description} 
	\item[Exponential Moment conditions:] 
The observation $Y$ satisfy to the following moment condition. There exists some positive numbers $v \geq \bE (Y_j^2)$ and $c$  such that
\begin{equation}\label{mc}
\forall j\geq 1,\: \forall k\geq 2,\quad \bE (Y_j^k)< \frac{k!}{2}vc^{k-2}.\end{equation}
\end{description}
These two conditions are required in order to obtain concentration bounds using first Hoeffding type inequality, then Bernstein inequality, see for instance \cite{geerb}. Requiring bounded SVD functions may be seen as a restrictive condition. Yet it is met when the eigenvectors are trigonometric functions. However, this condition can be also be turned into a moment condition if we replace the concentration bound by a Bernstein type inequality. Note also that the moment conditions on $Y$ amounts to require a bounded regression function $\vphi$ and equivalent moment conditions on the errors $U_j$. 
\begin{description}
	\item[IP: Degree of ill-posedness]
	We assume that there exists $t$, called the degree of ill-posedness of the operator which controls the decay of the eigenvalues of the operator $T$. More precisely, there are constants $\la_L,\la_U$ such that
	\begin{equation} \label{IP}
\la_L k^{-t}\leq	\la_k \leq \la_U k^{-t},\: \forall k\geq 1
	\end{equation} 
\end{description}
In this paper, we only consider the case of mildly ill-posed inverse problems, i.e when the eigenvalues decay at a polynomial rate. This assumption, also required in \cite{MR2158113}, is needed when comparing the residual error of the estimator with the risk in order to obtain the oracle inequality.    
\begin{description}
	\item[Enough ill-posedness : ]
	Let $\sigma_j^2= {\rm Var}(Y\psi_j(W))$. We assume that there exist two positive constants $\sigma^2_L$ and $\sigma_U^2$ such that
	\begin{equation} \label{cond:var} \forall j \geq 1,\quad \sigma_L^2 \leq \sigma^2_j \leq \sigma^2_U. \end{equation}
\end{description}
Note that Condition~\eqref{mc} implies the upper bound of Condition~\eqref{cond:var}. The lower bound is similar to the  variance condition in Assumption 3.1 in \cite{xiaoreiss}. We we also point out that this condition is not needed when building an estimator for the regression function. However it turns necessary when obtaining the lower bound to get a minimax result, or when obtaining an oracle inequality.

\subsection{Oracle inequality}
First, let $R_0(m,\vphi)$ be the quadratic estimation risk for the naive estimator  $\hat{\vphi}_m^0$ \eqref{naif}, defined by
\begin{align*}
R_0(m,\varphi) &= \bE \| \hat{\vphi}_m^0-\vphi \|^2 \\
& = \sum_{k>m} \varphi_k^2 + \frac{1}{n} \sum_{k=1}^m \lambda_k^{-2} \sigma_k^2, \ \forall m\in \dN. \end{align*}
The best model would be obtained by choosing a minimizer of this quantity, namely 
\begin{equation} \label{eq:m0}
m_0 = {\rm arg}\min_m R_0(m,\vphi). \end{equation}
This risk depends on the unknown function $\vphi$ hence $m_0$ is the oracle. We aim at constructing an estimator of $R_0(m,\vphi)$ which, by minimization, could give rise to a convenient choice for $m$, i.e as close as possible to $m_0$. The first step would be  to replace $\vphi_k$ by their estimates ${\hat{\la}_k}^{-1} \hat{r}_k$ and take for estimator of $\sigma^2_k$, $\hat\sigma_k^2$, defined by
\begin{align*}
\hat \sigma_k^2 &=\frac{1}{n} \sum_{i=1}^n \left(Y_i \psi_k(W_i) - \frac{1}{n}\sum_{i=1}^n Y_i \psi_k(W_i)\right)^2 \\
&=\frac{1}{n} \sum_{i=1}^n \left(Y_i \psi_k(W_i) - \hat{r}_k \right)^2. \end{align*} 
This would lead us to consider the empirical risk for any $m\leq M$, the cut-off which warrants a good behaviour for the $\hat \la_j$'s 
$$ U_0(m,r,\lambda)= - \sum_{k=1}^m \hat\lambda_k^{-2} \hat r_k^2 +  \frac{c}{n} \sum_{k=1}^m \hat\lambda_k^{-2} \hat\sigma_k^2, \: \forall m\in \dN, $$ for a well chosen constant $c$. The corresponding random oracle within the range of models which are considered would be 
\begin{equation}
\label{eq:m1}
m_1={\rm arg}\min_{ m\leq M}R_0(m,\vphi). \end{equation}
Unfortunately, the correlation between the errors $V_i$ and the observations $Y_i$ prevents an estimator defined as a minimizer of $U_0(m,r,\la)$ to achieve the quadratic risk $R_0(m,\vphi)$. Indeed, we have to use a stronger penalty, leading to an extra error in the estimation that shall be discussed later in the paper. More precisely, $c$ in the penalty is not a constant anymore but is allowed to depend on the number of observations $n$.\vskip .1in
Hence, now define $R(m,\vphi)$ the penalized estimation risk as
\begin{equation}
R(m,\varphi)  = \sum_{k>m} \varphi_k^2 + \frac{\log^2 n}{n} \sum_{k=1}^m \lambda_k^{-2} \sigma_k^2, \ \forall m\in \dN. \end{equation}
The best choice for $m$ would be a minimizer of this quantity, which yet depends on the unknown regression function $\vphi$. Hence, to mimic this risk, define the following empirical criterion
\begin{equation}
U(m,r,\lambda)= - \sum_{k=1}^m \hat\lambda_k^{-2} \hat r_k^2 +  \frac{\log^2 n}{n} \sum_{k=1}^m \hat\lambda_k^{-2} \hat\sigma_k^2, \: \forall m\in \dN.
\label{U(m,r,lambda)}
\end{equation}
Then, the best estimator is selected by minimizing this quantity as follows\begin{equation}
m^{\star} := \mathrm{arg} \min_{m \leq M} U(m,r,\lambda),
\label{mstar}
\end{equation}
Finally, the corresponding adaptive estimator $\varphi^{\star}$ is defined as:
\begin{equation}
\varphi^{\star} = \sum_{k=1}^{m^{\star}} \hat\lambda_k^{-1} \hat r_k \phi_k.
\label{varphistar}
\end{equation}
The performances of $\varphi^{\star}$ are presented in the following theorem.

\begin{thm}
Let $\varphi^{\star}$ the projection estimator defined in (\ref{varphistar}). Then, there exists $B_0,B_1,B_2$ and $\tau$ positive constants independent of $n$ such that:
\begin{eqnarray*}
\mathbb{E} \| \varphi^{\star} - \varphi \|^2 & \leq  & B_0 \log^2 (n). \left[ \inf_m R(m,\varphi) \right]+ \frac{B_1}{n} \left( \log (n). \|\varphi \|^2 \right)^{2 t} \\
& & \hspace{3cm} +\Omega + \log^2 (n) .\Gamma(\varphi),
\end{eqnarray*}
where $\Omega \leq B_2 (1+\| \varphi\|^2) \exp\left\lbrace -\log^{1+\tau} n \right\rbrace$, $m_0$ denotes the oracle bandwidth and
\begin{equation}
\Gamma(\varphi) = \sum_{k=\min(M_0,m_0)}^{m_0} \left[ \varphi_k^2 + \frac{1}{n} \lambda_k^{-2}\sigma_k^2 \right],
\label{gammaphi}
\end{equation} 
with the convention $\sum_a^b=0$ if $a=b$.
\end{thm}
We obtain a non asymptotic inequality which guarantees that the estimator achieves the optimal bound, up to a logarithmic factor, among all the estimators that could be constructed. We point out that we loss a $\log^2(n)$ factor when compared with the bound obtained in \cite{MR2158113}. The explanation of this loss comes from the fact that the error on the operator is not deterministic nor even due to a independent noisy observation of  the eigenvalues. Here, the $\la_k$'s have to be estimated using the available data by $\hat{\la}_k$. In the econometric model, both the operator and the regression function are estimated on the same sample, which leads to high correlation effects that are made explicit in Model \eqref{eq:modelreiss}, hampering the rate of convergence of the corresponding estimator.\\
\indent An oracle inequality only provides some information on the asymptotic behaviour of the estimator if the remainder term $\Gamma(\varphi)$ is of smaller order than the risk of the oracle. This remainder term models the error made when truncating the eigenvalues, i.e the error when selecting a model close to the random oracle $m_1 \leq M$ and not the true oracle $m_0$. In the next section, we prove that, under some assumptions, this extra term is smaller than the risk of the estimator. 
\subsection{Rate of convergence}
 To get a rate of convergence for the estimator, we need to specify the regularity of the unknown function $\vphi$ and compare it with the degree of ill-posedness of the operator $T$, following the usual conditions in the statistical literature on inverse problems, see for example \cite{Ruym96} or \cite{tsy}, \cite{MR2361904} for some examples. 
\begin{description}
	\item[Regularity Condition] 
	Assume that the function $\vphi$ is such that there exists $s$ and a constant $C$ such that
	\begin{equation}
	\sum_{k\geq 1} k^{2s}\vphi_k^2 < C
	\label{eq:Sobol}
\end{equation}
\end{description}
This Assumption corresponds to functions whose regularity is governed by the smoothness index $s$. This parameter is unknown and yet governs the rate of convergence. In the special cases where the eigenfunctions are the Fourier basis, this set corresponds to Sobolev classes. We prove that our estimator achieves the optimal rate of convergence without prior assumption on $s$.
\begin{cor} \label{larate}
Let $\vphi^\star$ be the model selection estimator defined in \eqref{varphistar}. Then, under the Sobolev embedding assumption \eqref{eq:Sobol}, we get the following rate of convergence
$$\mathbb{E} \| \varphi^{\star} - \varphi \|^2 = O\left(  \left(\frac{n}{ \log^{2\gamma} n}\right)^{\frac{-2s}{2s+2t+1}} \right), $$
with $\gamma=2+2s+2t$.\end{cor}
We point out that $\vphi^\star$ is constructed without prior knowledge of the unknown regularity $s$ of $\vphi$, yet achieving the optimal rate of convergence, up to some logarithmic terms. In this sense, our estimator is said to be asymptotically adaptive.
\begin{rem}
In an equivalent way, we could have imposed a supersmooth assumption, on the function $\vphi$, i.e assuming that for given $\gamma$, $t$ and constant $C$, $$\sum_{k=1}^\infty \exp(2 \gamma k^t) \vphi_k^2 < C.$$ Following the proof of Corollary~\ref{larate}, we obtain that $M_0 > m_0 \sim  (a 2 \gamma \log n)^{1/t}$ with $2 a \gamma >1$, leading to the optimal recovery rate for supersmooth functions in inverse problems.
\end{rem}
In conclusion, this work shows that provided the eigenvectors are known, for smooth functions $\vphi$, estimating the eigenvalues and using a threshold suffices to get a good estimator of the regression function in the instrumental variable framework. The price to pay for not knowing the operator is only an extra $\log^2 n$ with respect to usual inverse problems and is only due to the correlation induced by the $V_i$'s. One could object that when dealing with unknown operators, the knowledge of the eigenvectors is a huge hint and some papers have considered the case of completely unknown operators, using functional approach, see for instance \cite{florens}, \cite{joflo}, but their estimate clearly rely on smoothness assumptions for the regression. Hence the two approaches are complementary since we provide more refined adaptive result with the sake of stronger assumptions. Nevertheless, using similar techniques to develop a fully adaptive estimation procedure would be the last step towards a full understanding of the IV regression model.

\section{Technical lemmas}\label{s:appendix}
First of all, we point out that, throughout all the paper, $C$ denotes some generic constant that may vary from line to line.\vskip .1in
Recall that we have introduced
\begin{equation*}
M= \inf \left\lbrace k\leq N: |\hat\lambda_k | \leq \frac{1}{\sqrt{n}} \log n \right\rbrace -1,
\end{equation*}
The term $N$ provides a deterministic upper bound for $M$ and ensures that $M$ is not too large. Typically, choose $N=n^4$. The following lemma provides a control of the bandwidth $M$ by $M_0$ and $M_1$ respectively defined in (\ref{defM_0}) and (\ref{defM_1}).

\begin{lem}
\label{controlM}
Set $\mathcal{M}= \lbrace M_0 \leq M < M_1 \rbrace$. Then, for all $n\geq 1$,
$$ P(\mathcal{M}^c) \leq C M_0 e^{-\log^{1+\tau} n},$$
where $C$ and $\tau$ denote positive constants independent of $n$.
\end{lem}
\noindent
PROOF. It is easy to see that:
$$ P(\mathcal{M}^c) = P\left( \lbrace M<M_0 \rbrace \cup \lbrace M \geq M_1 \rbrace \right) \leq P(M<M_0) + P(M \geq M_1).$$
Using (\ref{M}) and (\ref{defM_1}),
$$ P(M \geq M_1)  =  P \left( \bigcap_{k=1}^{M_1} \left\lbrace |\hat\lambda_k | \geq \frac{1}{\sqrt{n}} \log n \right\rbrace \right)
         \leq  P \left(  |\hat\lambda_{M_1} | \geq \frac{1}{\sqrt{n}} \log n  \right).$$
Thanks to the definition of $\hat\lambda_{M_1}$: 
\begin{eqnarray*}
P(M \geq M_1)
& \leq & P \left(  \left|\hat\lambda_{M_1}- \lambda_{M_1} + \lambda_{M_1} \right| \geq                               \frac{1}{\sqrt{n}} \log n \right),\\ 
& \leq & P \left(  \left| \hat\lambda_{M_1}- \lambda_{M_1} \right| \geq   \frac{1}{\sqrt{n}} \log n                - |\lambda_{M_1} | \right) ,\\
& \leq &  P \left(  \left| \frac{1}{n} \sum_{i=1}^n  \phi_{M_1}(X_i)\psi_{M_1}(W_i) - \bE [\phi_{M_1}(X) \psi_{M_1}(W)] \right| \geq   b_n \right),
\end{eqnarray*}
where $b_n=n^{-1/2} \log n - |\lambda_{M_1}|$ for all $n\in \mathbb{N}$. Let $k\in\mathbb{N}$ and $x\in[0,1]$ be fixed. Assumption \eqref{basis} and Hoeffding inequality yields
\begin{eqnarray*}
P(| \hat\lambda_k - \lambda_k| > x) & \leq & 2 \exp \left\lbrace - \frac{(nx)^2}{2 \sum_{i=1}^n \mathrm{Var}(\phi_{M_1}(X_i)\psi_{M_1}(W_i))+ 2 nCx/3 } \right\rbrace,\\
	 & = & 2 \exp \left\lbrace - \frac{n x^2}{2 \mathrm{Var}(\phi_{M_1}(X)\psi_{M_1}(W))+ 2C x/3 } \right\rbrace.
\end{eqnarray*}
Using again the assumption \eqref{basis} on the bases $(\phi_k)_{k\in\dN}$ and $(\psi_k)_{k\in\dN}$,
$$ \mathrm{Var} ( \phi_{M_1}(X) \psi_{M_1}(W) ) \leq C_1^4 \mathbb{E} [ \phi_{M_1}^2(X) \psi_{M_1}^2(W) ] \leq 1.$$
Hence,
\begin{equation}
P(| \hat\lambda_k - \lambda_k | > x) \leq 2 \exp\left( - C \frac{x^2}{3} \right), \ \forall x\in [0,1],
\label{vpdeviation}
\end{equation}
with $C$ depending on $C_1$.\\
\indent Using (\ref{defM_1}), $1>b_n >0$ for all $n\in\dN$. Therefore, using (\ref{vpdeviation}) with $x=b_n$, we obtain:
\begin{eqnarray*}
P(M\geq M_1) \leq 2 \exp \left\lbrace -\frac{n b_n^2}{3} \right\rbrace & \leq & 2 \exp \left\lbrace -\frac{1}{3}(\log n - \log^{3/4} n )^2 \right\rbrace, \\
& \leq & C \exp \left\lbrace - \log^{1+\tau} n \right\rbrace,
\end{eqnarray*}
where $C$ and $\tau$ denote positive constants independent of $n$.\\
\\
The bound of $P(M < M_0)$ follows the same lines:
\begin{eqnarray*}
P(M < M_0) = P\left( \bigcup_{j=1}^{M_0} \left\lbrace |\hat\lambda_j | \leq \frac{\log n}{\sqrt{n}} \right\rbrace \right)  & \leq & \sum_{j=1}^{M_0}  P\left(  |\hat\lambda_j | \leq \frac{\log n}{\sqrt{n}} \right),\\
 & \leq & \sum_{j=1}^{M_0}  P\left(  \hat\lambda_j \leq \frac{\log n}{\sqrt{n}} \right).
\end{eqnarray*}
Let $j\in \lbrace 1,\dots, M_0 \rbrace$ be fixed.
\begin{eqnarray*}
P \left(  \hat\lambda_j \leq \frac{\log n}{\sqrt{n}} \right) & = & P\left(  \hat\lambda_j - \lambda_j \leq \frac{\log n}{\sqrt{n}} -\lambda_j \right), \\
& = & P\left(  \frac{1}{n} \sum_{i=1}^n \lbrace \phi_j(X_i)\psi_j(X_i) - \mathbb{E}[\phi_j(X_i)\psi_j(X_i)] \rbrace  \leq \tilde b_n \right),
\end{eqnarray*}
where $\tilde b_n= n^{-1/2}\log n - \lambda_j$ for all $n\in\dN$. Thanks to (\ref{defM_0}), $\tilde b_n <0$ for all $n\in\dN$. Using Hoeffding inequality and Assumption \eqref{basis} :
$$ P \left(  \hat\lambda_j \leq \frac{\log n}{\sqrt{n}} \right) \leq \exp \left\lbrace - \frac{n \tilde b_n^2}{2 + 2/3 |\tilde b_n |} \right\rbrace \leq C \exp \left\lbrace - \log^{1+\tau} n \right\rbrace,$$
for some $C,\tau>0$. This concludes the proof of Lemma \ref{controlM}.
\begin{flushright}
$\Box$
\end{flushright}

\begin{lem}
\label{controlB}
Let $\mathcal{B}$ the event defined by:
$$ \mathcal{B} = \bigcap_{k=1}^M \left\lbrace | \lambda_k^{-1} \mu_k | \leq \frac{1}{2} \right\rbrace, \ \mathrm{where} \ \mu_k= \hat\lambda_k - \lambda_k, \ \forall k\in \dN^*.$$
Then,
$$P(\mathcal{B}^c) \leq CM_1 e^{-\log^{1+\tau} n},$$
for some $\tau>0$ and positive constant $C$.
\end{lem}
\noindent
PROOF. Using simple algebra and Lemma \ref{controlM}
\begin{eqnarray*}
P(\mathcal{B}^c) & = & P(\mathcal{B}^c \cap \mathcal{M})+P(\mathcal{B}^c \cap \mathcal{M}^c),\\
                 & \leq & P(\mathcal{B}^c \cap \mathcal{M})+P(\mathcal{M}^c),\\
		 & \leq & P(\mathcal{B}^c \cap \mathcal{M})+C M_0 e^{-\log^{1+\tau} n}.
\end{eqnarray*}
Then,
$$ P(\mathcal{B}^c \cap \mathcal{M})  =  P\left( \bigcup_{k=1}^M \left\lbrace | \lambda_k^{-1}\mu_k | > \frac{1}{2} \right\rbrace \cap \mathcal{M} \right) \leq  P\left( \bigcup_{k=1}^{M_1-1} \left\lbrace | \lambda_k^{-1}\mu_k | \geq \frac{1}{2} \right\rbrace\right).$$
Let $k\in \lbrace 1,\dots, M_1-1 \rbrace$ be fixed. Remark that:
$$ P\left( | \lambda_k^{-1}\mu_k | \geq \frac{1}{2} \right)  =  P\left( | \mu_k | \geq \frac{|\lambda_k |}{2} \right)
   \leq  P \left( |\hat\lambda_k - \lambda_k | \geq \frac{1}{2\sqrt{n}} \log^{3/4} n \right).$$
Then, using (\ref{vpdeviation}) with $x=2n^{-1/2} \log^{3/4} n$:
\begin{equation}
P \left( |\hat\lambda_k - \lambda_k | \geq \frac{1}{2\sqrt{n}} \log^{3/4} n \right) \leq C e^{-\log^{1+\tau} n},
\label{bornemu}
\end{equation}
for some $\tau>0$ and a positive constant $C$. This concludes the proof of Lemma \ref{controlB}.
\begin{flushright}
$\Box$
\end{flushright}

The following lemma provides some tools for the control of the ratio $\hat\lambda_k^{-1} \lambda_k$ on the event $\mathcal{B}$.

\begin{lem}
\label{ratiovp}
For all $k\leq M$, we have:
$$ \left( \frac{\lambda_k}{\hat\lambda_k} -1 \right)^2 \mathbf{1}_{\mathcal{B}} \leq \frac{2}{3}\lambda_k^{-2} (\hat\lambda_k - \lambda_k)^2 \mathbf{1}_{\mathcal{B}}.$$
Moreover, we have the following expansion:
$$\frac{\hat\lambda_k^{-1}}{\lambda_k} = 1 - \lambda_k^{-1}(\hat\lambda_k-\lambda_k) + \lambda_k^{-2} (\hat\lambda_k - \lambda_k)^2 \nu_k,$$
where $\nu_k$ is uniformly bounded on the event $\mathcal{B}$.
\end{lem}
\noindent
PROOF. Let $k\leq M$ be fixed. Then
$$ \left( \frac{\lambda_k}{\hat\lambda_k} -1 \right)^2 \mathbf{1}_{\mathcal{B}} = \left( \frac{\mu_k}{\hat\lambda_k} \right)^2 \mathbf{1}_{\mathcal{B}} = \left( \frac{\mu_k}{\lambda_k+ \mu_k} \right)^2 \mathbf{1}_{\mathcal{B}} \leq \frac{2}{3}\lambda_k^{-2} (\hat\lambda_k - \lambda_k)^2 \mathbf{1}_{\mathcal{B}},$$
where the $\mu_k$ are defined in Lemma \ref{controlB}. The end of the proof is based on a Taylor expansion of the ratio $\hat\lambda_k^{-1} \lambda_k= (1+\lambda_k^{-1}\mu_k)^{-1}$. The variable $\nu_k$ depends on $\lambda_k^{-1}\mu_k$ and can be easily bounded on the event $\mathcal{B}$.
\begin{flushright}
$\Box$
\end{flushright}

\begin{lem}
\label{T1}
Let $\bar m$ a random variable measurable with respect to $(Y_i,X_i,W_i)_{i=1\dots n}$ such that $\bar m \leq M$. Then, for all $K>1$ and $\gamma>0$,
$$ (i) \ \bE \left[\sum_{k=1}^{\bar m} \hat\lambda_k^{-2} (\hat r_k -r_k)^2 \right] \leq \frac{\log^K (n)}{n} \bE \left[\sum_{k=1}^{\bar m} \hat\lambda_k^{-2} \sigma_k^2 \right]+ CNn e^{-\log^K n},$$
\begin{eqnarray*}
(ii) \ \bE \left[\sum_{k=1}^{\bar m} \lambda_k^{-2} (\hat r_k -r_k)r_k \right] 
& \leq & \gamma^{-1} \frac{\log^K (n)}{n} \bE \left[\sum_{k=1}^{\bar m} \hat\lambda_k^{-2} \sigma_k^2 \right]+ C\gamma^{-1}N^{2t+1} e^{-\log^K n}\\
& & + \gamma^{-1} R(m_0,\varphi) + \gamma \bE \sum_{k> \bar m} \varphi_k^2,
\end{eqnarray*}
where $C>0$ is a positive constant independent of $n$, $m_0$ denotes the oracle bandwidth and $N$ has been introduced in (\ref{M}).
\end{lem}

\noindent
PROOF. Let $Q>0$ a positive term which will be chosen later. With simple algebra:
\begin{eqnarray}
\lefteqn{\bE \left[\sum_{k=1}^{\bar m} \hat\lambda_k^{-2} (\hat r_k -r_k)^2 \right]} \nonumber \\
& = & \bE \sum_{k=1}^{\bar m} \hat\lambda_k^{-2} (\hat r_k -r_k)^2 \mathbf{1}_{\left\lbrace (\hat r_k -r_k)^2 < \frac{Q\sigma_k^2}{n} \right\rbrace} + \bE \sum_{k=1}^{\bar m} \hat\lambda_k^{-2} (\hat r_k -r_k)^2 \mathbf{1}_{\left\lbrace (\hat r_k -r_k)^2 \geq \frac{Q\sigma_k^2}{n} \right\rbrace}, \nonumber \\
& \leq & \frac{Q}{n} \bE \left[\sum_{k=1}^{\bar m} \hat\lambda_k^{-2} \sigma_k^2 \right] + \bE \sum_{k=1}^{\bar m} \hat\lambda_k^{-2} (\hat r_k -r_k)^2 \mathbf{1}_{\left\lbrace (\hat r_k -r_k)^2 \geq \frac{Q\sigma_k^2}{n} \right\rbrace}.
\label{l1}
\end{eqnarray}
In the sequel, we are interested in the behavior of the second term in the right hand side of (\ref{l1}). Since $\hat\lambda_k^{-2} \leq n \log^{-2} n$ for all $k\leq M$ and $\bar m \leq N$, we obtain:
\begin{equation}
\bE \sum_{k=1}^{\bar m} \hat\lambda_k^{-2} (\hat r_k -r_k)^2 \mathbf{1}_{\left\lbrace (\hat r_k -r_k)^2 \geq \frac{Q\sigma_k^2}{n} \right\rbrace} \leq \frac{n}{\log^2 n} \sum_{k=1}^N \bE (\hat r_k -r_k)^2 \mathbf{1}_{\left\lbrace (\hat r_k -r_k)^2 \geq \frac{Q\sigma_k^2}{n} \right\rbrace}.
\label{bornesomme}
\end{equation}
Let $k\in \lbrace 1,\dots,N \rbrace$ be fixed. It follows from integration by part that:
$$ \bE (\hat r_k -r_k)^2 \mathbf{1}_{\left\lbrace (\hat r_k -r_k)^2 \geq \frac{Q\sigma_k^2}{n} \right\rbrace} \leq \int_{\frac{Q\sigma_k^2}{n}}^{+\infty} P\left( (\hat r_k -r_k)^2 >x \right) dx.$$
Then,
$$P\left( (\hat r_k-r_k)^2 \geq x \right) = P\left( |\hat r_k - r_k | \geq \sqrt{x} \right).$$
Assumption \eqref{mc} together with Bernstein inequality entails that:
\begin{eqnarray*}
P\left( |\hat r_k - r_k | \geq \sqrt{x} \right)
  & = & P\left( \left| \frac{1}{n} \sum_{i=1}^n (Y_i\psi_k(W_i)- \bE[Y_i\psi_k(W_i)]) \right| \geq \sqrt{x} \right),\\
  & \leq & \exp \left\lbrace - \frac{n^2x}{2\sum_{i=1}^n \mathrm{Var}(Y_i \psi_k(W_i)) + Cn \sqrt{x}} \right\rbrace,\\
  & = & \exp \left\lbrace - \frac{nx}{2\sigma_k^2 + C \sqrt{x}} \right\rbrace.
\end{eqnarray*}
Now remark that:
$$ 2 \sigma_k^2 = C\sqrt{x} \Leftrightarrow x = D, \ \mathrm{with} \ D=\left( \frac{2\sigma_k^2}{C} \right)^2.$$
We obtain:
\begin{eqnarray*}
\lefteqn{\bE (\hat r_k -r_k)^2 \mathbf{1}_{\left\lbrace (\hat r_k -r_k)^2 \geq \frac{Q\sigma_k^2}{n} \right\rbrace}} \\
& \leq & \int_{\frac{Q\sigma_k^2}{n}}^D \exp \left\lbrace - \frac{nx}{2\sigma_k^2 + C \sqrt{x}} \right\rbrace dx + \int_D^{+\infty} \exp \left\lbrace - \frac{nx}{2\sigma_k^2 + C \sqrt{x}} \right\rbrace dx, \\
& \leq & \int_{\frac{Q\sigma_k^2}{n}}^D \exp \left\lbrace - \frac{nx}{4\sigma_k^2} \right\rbrace dx + \int_D^{+\infty} \exp \left\lbrace - \frac{nx}{C \sqrt{x}} \right\rbrace dx ,\\
& \leq & \left[ - \frac{4\sigma_k^2}{n} e^{-\frac{nx}{4\sigma_k^2}} \right]_{Q\sigma_k^2/n}^{+\infty} + \int_D^{+\infty} \exp \left\lbrace - Cn\sqrt{x} \right\rbrace dx, \\
& \leq & \frac{4\sigma_k^2}{n} \exp\left\lbrace - \frac{n}{4\sigma_k^2} \frac{Q\sigma_k^2}{n} \right\rbrace + \frac{\sqrt{D}+1}{Cn} e^{-C\sqrt{D} n} ,\\
& \leq & \frac{4\sigma_k^2}{n} e^{-Q/4} + e^{-Cn}.
\end{eqnarray*}
Hence, we have
\begin{equation}
\bE (\hat r_k -r_k)^2 \mathbf{1}_{\left\lbrace (\hat r_k -r_k)^2 \geq \frac{Q\sigma_k^2}{n} \right\rbrace} \leq \frac{C\sigma_k^2}{n} e^{-Q/4} + e^{-Cn},
\label{borneT_1+T_2}
\end{equation}
for some $C>0$. Using (\ref{borneT_1+T_2}) and (\ref{bornesomme}),
$$\bE \sum_{k=1}^{\bar m} \hat\lambda_k^{-2} (\hat r_k -r_k)^2 \mathbf{1}_{\left\lbrace (\hat r_k -r_k)^2 \geq \frac{Q\sigma_k^2}{n} \right\rbrace} \leq \frac{C Nn}{\log^2 n} e^{-Q/4}+ e^{-Cn}.$$
From (\ref{l1}), we eventually obtain:
$$\bE \left[\sum_{k=1}^{\bar m} \hat\lambda_k^{-2} (\hat r_k -r_k)^2 \right]
  \leq  \frac{Q}{n} \bE \left[\sum_{k=1}^{\bar m} \hat\lambda_k^{-2} \sigma_k^2 \right] + \frac{C Nn}{\log^2 n} e^{-Q/4} + e^{-Cn}.$$
Choose $Q=\log^K (n)$ in order to conclude the proof of $(i)$.\\
\\
 Now, consider the bound of $(ii)$. Let $m_0$ the oracle bandwidth defined in (\ref{eq:m0}). With the convention $\sum_a^b=-\sum_b^a$ if $b<a$,
\begin{eqnarray}
\bE \sum_{k=1}^{\bar m} \lambda_k^{-2} (\hat r_k-r_k)r_k & = & \bE \sum_{k=m_0}^{\bar m} \lambda_k^{-2} (\hat r_k-r_k)r_k, \nonumber \\
    & \leq & \bE \left| \sum_{k=m_0}^{\bar m} \lambda_k^{-2} (\hat r_k-r_k)r_k \right| ,\nonumber \\
    & \leq & \bE \sum_{k=1}^{+\infty} \left| (\mathbf{1}_{\lbrace k\leq \bar m \rbrace} - \mathbf{1}_{\lbrace k\leq m_0 \rbrace}) \lambda_k^{-2} (\hat r_k-r_k)r_k \right|.
\label{indicatrices}
\end{eqnarray}
Indeed, $\bE [ \hat r_k] = r_k$ for all $k\in \mathbb{N}$. Then remark that:
\begin{eqnarray}
\left| \mathbf{1}_{\lbrace k\leq \bar m \rbrace} - \mathbf{1}_{\lbrace k\leq m_0 \rbrace} \right| 
& = & \left| (\mathbf{1}_{\lbrace k\leq \bar m \rbrace} + \mathbf{1}_{\lbrace k\leq m_0 \rbrace}) (\mathbf{1}_{\lbrace k\leq \bar m \rbrace} - \mathbf{1}_{\lbrace k\leq m_0 \rbrace}) \right|, \nonumber \\
& = & (\mathbf{1}_{\lbrace k\leq \bar m \rbrace} + \mathbf{1}_{\lbrace k\leq m_0 \rbrace}) \left| \mathbf{1}_{\lbrace k > \bar m \rbrace} - \mathbf{1}_{\lbrace k> m_0 \rbrace} \right|, \nonumber \\
& \leq & \mathbf{1}_{\lbrace k> \bar m \rbrace} \mathbf{1}_{\lbrace k\leq m_0 \rbrace} + \mathbf{1}_{\lbrace k> m_0 \rbrace} \mathbf{1}_{\lbrace k\leq \bar m \rbrace}.
\label{indicatrices2}
\end{eqnarray}
Using the Cauchy-Schwartz inequality and using that for all  $a,b$ and $1>\gamma>0$, $2ab\leq \gamma a^2 +\gamma^{-1} b^2$:
\begin{eqnarray*}
\lefteqn{\bE \sum_{k=1}^{\bar m} \lambda_k^{-2} (\hat r_k-r_k)r_k} \\ 
  & \leq &  \sqrt{ \bE\sum_{k> \bar m} \lambda_k^{-2} r_k^2} \sqrt{\bE\sum_{k\leq m_0} \lambda_k^{-2} (\hat r_k - r_k)^2}  +  \sqrt{ \bE\sum_{k>m_0} \lambda_k^{-2} r_k^2} \sqrt{\bE\sum_{k\leq \bar m} \lambda_k^{-2} (\hat r_k - r_k)^2}, \\
  & \leq & \gamma \left\lbrace \bE \sum_{k> \bar m} \varphi_k^2 + \sum_{k>m_0} \varphi_k^2 \right\rbrace + \gamma^{-1} \left\lbrace \bE \sum_{k=1}^{\bar m} \lambda_k^{-2} (\hat r_k - r_k)^2 + \bE \sum_{k=1}^{m_0} \lambda_k^{-2} (\hat r_k - r_k)^2 \right\rbrace.
\end{eqnarray*}
We eventually obtain:
$$\bE \sum_{k=1}^{\bar m} \lambda_k^{-2} (\hat r_k-r_k)r_k
\leq  \gamma^{-1} R(m_0,\varphi) + \gamma \bE \sum_{k> \bar m} \varphi_k^2 + \gamma^{-1} \left\lbrace \bE \sum_{k=1}^{\bar m} \lambda_k^{-2} (\hat r_k - r_k)^2 \right\rbrace.$$
We conclude the proof using a similar to $(i)$ string of inequalities. In particular, using Assumption (\ref{IP}), we obtain the bound $\lambda_k^{-2} \leq CN^{2t}$ for all $k\leq M$.
\begin{flushright}
$\Box$
\end{flushright}

\begin{lem}
Let $\bar m$ a random variable measurable with respect to $(Y_i,X_i,W_i)_{i=1\dots n}$ such that $\bar m \leq M$. Then, for all $\gamma \in (0,1)$,
\begin{eqnarray*}
\bE \sum_{k=1}^{\bar m} ( \hat\lambda_k^{-2} - \lambda_k^{-2} ) r_k^2 & \leq & \frac{\gamma + \gamma^{-1} \log^{3/2} n}{n} \bE \left[ \sum_{k=1}^{\bar m} \lambda_k^{-2} \sigma_k^2 \right]+ \frac{1}{n}\left( \frac{\log^2 (n). \|\varphi \|^2}{\gamma} \right)^{2t} \\
& & \hspace{2cm} + \log^2(n). R(m_0,\varphi)+ \Omega.
\end{eqnarray*}
\label{risquevpbruitee}
\end{lem}
\noindent
PROOF. The term in the left hand side can be rewritten as:
$$ \bE \sum_{k=1}^{\bar m} ( \hat\lambda_k^{-2} - \lambda_k^{-2} ) r_k^2 = \bE \sum_{k=1}^{\bar m} \left( \frac{\lambda_k^2}{\hat\lambda_k^2} - 1 \right) \lambda_k^{-2} r_k^2 = \bE \sum_{k=1}^{\bar m} \left( \frac{\lambda_k^2}{\hat \lambda_k^2} -1 \right) \varphi_k^2.$$
Using Lemma \ref{ratiovp}, we obtain:
$$ \bE \sum_{k=1}^{\bar m} ( \hat\lambda_k^{-2} - \lambda_k^{-2} ) r_k^2
=  - \bE \left[ \sum_{k=1}^{\bar m} \varphi_k^2 \lambda_k^{-1}\mu_k \right] + \bE \left[ \sum_{k=1}^{\bar m} \varphi_k^2 \lambda_k^{-2}\mu_k^2 \nu_k \right] =  W_1 + W_2,$$
where the $\mu_k$ are defined in Lemma \ref{controlB}. First consider the bound of $W_2$. Using (\ref{vpdeviation}) with $x=n^{-1/2} \log n$, we obtain:
\begin{eqnarray}
W_2 & = & \bE \left[ \sum_{k=1}^{\bar m} \varphi_k^2 \lambda_k^{-2} \mu_k^2 \nu_k \right]
 \leq C \bE \left[ \sum_{k=1}^{\bar m} \varphi_k^2 \lambda_k^{-2} \mu_k^2 \right] + \Omega, \nonumber \\
& \leq & C\frac{\log^2 n}{n} \bE \left[ \sum_{k=1}^{\bar m} \varphi_k^2 \lambda_k^{-2} \right] + C \|\varphi\|^2 e^{-\log^{1+\tau} n},
\label{boundW2}
\end{eqnarray}
where $C$ denotes a positive constant independent of $n$. Thanks to our assumptions on the sequence $(\lambda_k)_{k\in\mathbb{N}}$, for all $\gamma>0$
\begin{eqnarray}
W_2 & \leq & \frac{\log^2 n}{n} \|\varphi\|^2 \bE \sup_{k\leq \bar m} \lambda_k^{-2} + C \|\varphi\|^2 e^{-\log^{1+\tau} n}, \nonumber\\
& \leq & \frac{\gamma}{n} \sum_{k=1}^{\bar m} \lambda_k^{-2}\sigma_k^2 + \frac{C}{n}\left( \frac{\log^2 (n). \|\varphi \|^2}{\gamma} \right)^{2t} + \Omega,
\label{W2}
\end{eqnarray}
where for the last inequality, we have used (\ref{IP}) and the bound:
$$ \sup_{k\leq \bar m} \lambda_k^{-2} \leq \frac{1}{x} \sum_{k=1}^{\bar m} \lambda_k^{-2} + Cx^{2t},$$
with $x= \gamma^{-1} \log^2(n). \|\varphi\|^2$. More details on this bound can be found in \cite{cgpt}.\\

We are now interested in the bound of $W_1$. Using (\ref{indicatrices2}) and a similar to (\ref{indicatrices}) string of inequalities, we obtain:
\begin{eqnarray*}
W_1 & = & \bE \sum_{k=1}^{\bar m} \varphi_k^2 \lambda_k^{-2} \mu_k,\\
    & \leq & \bE \sum_{k=1}^{+\infty} \mathbf{1}_{\lbrace k> \bar m \rbrace} \mathbf{1}_{\lbrace k\leq m_0 \rbrace} \varphi_k^2 |\lambda_k^{-1} \mu_k | + \bE \sum_{k=1}^{+\infty} \mathbf{1}_{\lbrace k> m_0 \rbrace} \mathbf{1}_{\lbrace k\leq \bar m \rbrace} \varphi_k^2 |\lambda_k^{-1} \mu_k |,\\
& \leq & \sqrt{ \bE \sum_{k> \bar m} \varphi_k^2} \sqrt{ \bE \sum_{k\leq m_0} \lambda_k^{-2} (\hat \lambda_k - \lambda_k)^2} + \sqrt{ \bE \sum_{k> m_0} \varphi_k^2} \sqrt{ \bE \sum_{k\leq \bar m} \lambda_k^{-2} (\hat \lambda_k - \lambda_k)^2}.
\end{eqnarray*}
Hence, for all $\gamma> 0$,
$$ W_1 \leq \gamma \left\lbrace \bE \sum_{k > \bar m} \varphi_k^2 + \sum_{k > m_0} \varphi_k^2 \right\rbrace + \gamma^{-1} \left\lbrace \bE \sum_{k=1}^{\bar m} \lambda_k^{-2} (\hat \lambda_k -  \lambda_k)^2 + \bE \sum_{k=1}^{m_0} \lambda_k^{-2} (\hat \lambda_k - \lambda_k)^2 \right\rbrace.$$ 
Using (\ref{vpdeviation}) once again with $x=n^{-1/2} \log^{3/4} n$, we obtain for all $\gamma>0$:
$$ W_1 \leq \gamma \left\lbrace \bE \sum_{k > \bar m} \varphi_k^2 + \sum_{k > m_0} \varphi_k^2 \right\rbrace + \frac{\gamma^{-1} \log^2 n}{n} \left\lbrace \bE \sum_{k=1}^{\bar m} \lambda_k^{-2} \sigma_k^2 + \sum_{k=1}^{m_0} \lambda_k^{-2} \sigma_k^2 \right\rbrace.$$ 
This concludes the proof of Lemma 3.5.
\begin{flushright}
$\Box$
\end{flushright}

\begin{lem}
Let $\bar m$ a random variable measurable with respect to $(Y_i,X_i,W_i)_{i=1\dots n}$ such that $\bar m \leq M$. Then,
$$ \frac{1}{n} \bE \left[ \sum_{k=1}^{\bar m} \hat\lambda_k^{-2} (\hat \sigma_k^2 - \sigma_k^2) \right] \leq C\frac{\log n}{n^{3/2}}. \bE \left[ \sum_{k=1}^{\bar m} \hat\lambda_k^{-2} \sigma_k^2 \right] + \frac{1}{n} \bE \left[ \sum_{k=1}^{\bar m} \hat \lambda_k^{-2} (r_k^2 - \hat r_k^2) \right]+ C e^{-\log^2 n},$$
for some $C>0$ independent of $n$.
\end{lem}
\noindent
PROOF. First remark that, for all $k\geq 1$,
\begin{eqnarray*}
\hat\sigma_k^2 - \sigma_k^2 & = & \frac{1}{n} \sum_{i=1}^n (Y_i \psi_k(W_i) - \hat r_k)^2 - \sigma_k^2, \\
    & = & \frac{1}{n} \sum_{i=1}^n Y_i^2 \psi_k^2(W_i) + \hat r_k^2 - \frac{2\hat r_k}{n} \sum_{i=1}^n Y_i \psi_k(W_i) - \sigma_k^2,\\
    & = & \frac{1}{n} \sum_{i=1}^n Y_i^2 \psi_k^2(W_i) + \hat r_k^2 - 2\hat r_k^2 - \left( \bE[Y^2 \psi_k^2(W)] - \bE[Y \psi_k(W)]^2 \right),\\
    & = & \frac{1}{n} \sum_{i=1}^n \left\lbrace Y_i^2 \psi_k^2(W_i) - \bE[Y^2 \psi^2_k(W)] \right\rbrace + (r_k^2 - \hat r_k^2).
\end{eqnarray*}
Hence, we obtain
\begin{equation}
\frac{1}{n} \bE \left[ \sum_{k=1}^{\bar m} \hat\lambda_k^{-2} (\hat \sigma_k^2 - \sigma_k^2) \right] =  \frac{1}{n} \bE \left[ \sum_{k=1}^{\bar m} \hat\lambda_k^{-2} \rho_k \right] + \frac{1}{n} \bE \left[ \sum_{k=1}^{\bar m} \hat\lambda_k^{-2} ( r_k^2 - \hat r_k^2) \right] ,
\label{1lem4}
\end{equation}
where for all $k\in \mathbb{N}$:
$$\rho_k= \frac{1}{n} \sum_{i=1}^n \left\lbrace Y_i^2 \psi_k^2(W_i) - \bE[Y^2 \psi_k(W)] \right\rbrace.$$
We are interested in the first term in the right hand side of (\ref{1lem4}). Let $\delta>0$ a positive constant which will be chosen later:
\begin{eqnarray*}
\frac{1}{n} \bE \left[ \sum_{k=1}^{\bar m} \hat\lambda_k^{-2} \rho_k \right]
& = & \frac{1}{n} \bE \left[ \sum_{k=1}^{\bar m} \hat\lambda_k^{-2} \rho_k \mathbf{1}_{\lbrace \rho_k \leq \delta \rbrace} \right]+ \frac{1}{n} \bE \left[ \sum_{k=1}^{\bar m} \hat\lambda_k^{-2} \rho_k \mathbf{1}_{\lbrace \rho_k > \delta \rbrace} \right],\\
& \leq & \frac{\delta}{n} \bE \left[ \sum_{k=1}^{\bar m} \hat\lambda_k^{-2} \right] + \frac{1}{n} \bE \left[ \sum_{k=1}^{\bar m} \hat\lambda_k^{-2} \rho_k \mathbf{1}_{\lbrace \rho_k > \delta \rbrace} \right].
\end{eqnarray*}
Since $\bar m\leq M$, from integration by part,
$$\frac{1}{n} \bE \left[ \sum_{k=1}^{\bar m} \hat\lambda_k^{-2} \rho_k \mathbf{1}_{\lbrace \rho_k > \delta \rbrace} \right]
\leq  \frac{1}{\log^2 n} \sum_{k=1}^{N} \bE \rho_k \mathbf{1}_{\lbrace \rho_k > \delta \rbrace} = \frac{1}{\log^2 n}  \sum_{k=1}^{N} \int_{\delta}^{+\infty} P(\rho_k \geq x) dx.$$
Let $k\in \dN$ and $x\geq \delta$ be fixed. Using Bernstein inequality:
\begin{eqnarray*}
P(\rho_k \geq x) & = & P\left( \frac{1}{n} \sum_{i=1}^n \left\lbrace Y_i^2 \psi_k^2(W_i) - \bE[Y^2 \psi_k(W)] \right\rbrace \geq x \right), \\
   & \leq & \exp \left\lbrace - \frac{n^2x^2}{2 \sum_{i=1}^n \mathrm{Var}(Y_i^2 \psi_k^2(W_i)) + Cxn/3} \right\rbrace, \\
   & \leq & \exp \left\lbrace - \frac{n^2 x^2}{2nD_0 + D_1 nx} \right\rbrace,\\
   & \leq & \exp \left\lbrace - \frac{n x^2}{2D_0 + D_1 x} \right\rbrace,
\end{eqnarray*}
with the hypotheses \eqref{mc} and \eqref{basis} on $Y$ and $(\psi_k)_k$. The constants $D_0$ and $D_1$ are positive and independent of $n$. Therefore, for all $k\leq N$,
\begin{eqnarray*}
\lefteqn{\int_{\delta}^{+\infty} P(\rho_k \geq x) dx }\\
  & \leq & \int_{\delta}^{2D_0/D_1} \exp \left\lbrace - \frac{n x^2}{2D_0 + D_1 x} \right\rbrace dx + \int_{2D_0/D_1}^{+\infty} \exp \left\lbrace - \frac{n x^2}{2D_0 + D_1 x} \right\rbrace dx,\\
  & \leq & \int_{\delta}^{2D_0/D_1} \exp \lbrace - Cnx^2 \rbrace dx + \int_{2D_0/D_1}^{+\infty} \exp \lbrace - nx \rbrace dx,\\
  & \leq & \int_{\delta}^{+\infty} \exp \lbrace - Cn\delta x \rbrace dx + \frac{1}{n} e^{-Cn},\\
  & \leq & \frac{C}{n\delta} \exp \lbrace -n\delta^2 \rbrace + n^{-1} e^{-Cn},
\end{eqnarray*}
for some $C>0$. Choosing $\delta = n^{-1/2} \log n$ and using Assumption (\ref{cond:var}), we obtain:
$$ \frac{1}{n} \bE \left[ \sum_{k=1}^{\bar m} \hat\lambda_k^{-2} \rho_k \right] \leq C\frac{\log n}{n^{3/2}} \bE \left[ \sum_{k=1}^{\bar m} \hat\lambda_k^{-2} \sigma_k^2 \right] + C e^{-\log^2 n}.$$
We use (\ref{1lem4}) in order to conclude the proof.
\begin{flushright}
$\Box$
\end{flushright}

\section{Proofs}
\textsc{Proof of Theorem 1}. The proof of our main result can be decomposed in four steps. In a first time, we prove that the quadratic risk of $\varphi^{\star}$ is close, up to some residual terms, to $\bE \bar R(m^{\star},\varphi)$ where
\begin{equation}
\bar R(m,\varphi) = \sum_{k>m} \varphi_k^2 + \frac{\log^2 n}{n} \sum_{k=1}^m \hat\lambda_k^{-2} \sigma_k^2, \ \forall m\in \dN.
\label{R(m,varphi)}
\end{equation}
This result is uniform in $m$ and justifies our choice of $\bar R(m,\varphi)$ as a criterion for the bandwidth selection.

In a second time, we show that $\bE \bar R(m^{\star},\varphi)$ and $\bE U(m^{\star},r,\varphi)$ are in some sense comparable. Then, according to the definition of $m^{\star}$ in (\ref{mstar}),
$$ U(m^{\star},r,\varphi) \leq U(m,r,\varphi), \forall m\leq M.$$
We will conclude the proof by proving that for all $m\leq M$, $\bE U(m,r,\varphi) = \bE \| \hat\varphi_m - \varphi \|^2$, up to a log term and some residual terms.\\

In order to begin the proof, remark that:
$$ \bE \|\varphi^{\star} - \varphi \|^2  =  \bE \sum_{k=1}^{+\infty} (\varphi_k^{\star} - \varphi_k)^2
    =  \bE \sum_{k> m^{\star}} \varphi_k^2 + \bE \sum_{k=1}^{m^{\star}} (\hat\lambda_k^{-1} \hat r_k - \varphi_k)^2.$$
This is the usual bias-variance decomposition. Then
\begin{eqnarray*}
\bE \sum_{k=1}^{m^{\star}} (\hat\lambda_k^{-1} \hat r_k - \varphi_k)^2
& = & \bE \sum_{k=1}^{m^{\star}} (\hat\lambda_k^{-1} \hat r_k -\hat\lambda_k^{-1} r_k+ \hat\lambda_k^{-1}r_k - \varphi_k)^2,\\
& \leq & 2 \bE \sum_{k=1}^{m^{\star}} \hat\lambda_k^{-2} (\hat r_k -r_k)^2 + 2\bE \sum_{k=1}^{m^{\star}}   (\hat\lambda_k^{-1}r_k - \varphi_k)^2 = T_1 + T_2.
\end{eqnarray*}
Concerning $T_2$, we use the following approach. For all $\gamma>0$, using Lemma 3.3 and the bounds (\ref{boundW2}) and (\ref{W2}):
\begin{eqnarray}
T_2 & = & \bE \sum_{k=1}^{m^{\star}} (\hat\lambda_k^{-1} r_k - \varphi_k)^2
     =  \bE \sum_{k=1}^{m^{\star}} \left( \frac{\lambda_k}{\hat\lambda_k} -1 \right)^2 \varphi_k^2, \nonumber \\
    & = & \bE \sum_{k=1}^{m^{\star}} \left( \frac{\lambda_k}{\hat\lambda_k} -1 \right)^2 \varphi_k^2 \mathbf{1}_{\mathcal{B}} + \bE \sum_{k=1}^{m^{\star}} \left( \frac{\lambda_k}{\hat\lambda_k} -1 \right)^2 \varphi_k^2 \mathbf{1}_{\mathcal{B}^c}, \nonumber \\
    & \leq & \frac{2}{3} \bE \left[ \sum_{k=1}^{m^{\star}} \lambda_k^{-2} \mu_k^2 \varphi_k^2 \right] + \Omega, \nonumber \\
    & \leq & \frac{\gamma}{n} \bE \sum_{k=1}^{m^{\star}} \lambda_k^{-2}\sigma_k^2 + C \left( \frac{\|\varphi\|^2 \log^2 (n)}{\gamma} \right)^{2t} + \Omega.
\label{numT2} 
\end{eqnarray}
where $\mu_k= \hat\lambda_k - \lambda_k$ for all $k\in \mathbb{N}$. The term $T_1$ is bounded using Lemma \ref{T1} with $\bar m =m^{\star}$ and $K=2$. Hence, for all $\gamma>0$,
\begin{equation}
\bE \| \varphi^{\star} - \varphi\|^2 \leq (1+\gamma) \bE \bar R(m^{\star},\varphi) + \frac{C}{n} \left( \frac{\log^2 (n). \|\varphi \|^2}{\gamma} \right)^{2t} + \Omega,
\label{Etape1}
\end{equation}
where $\bar R(m^{\star},\varphi)$ is introduced in (\ref{R(m,varphi)}). This concludes the first step of our proof. \\

Now, our aim is to write $\bE \bar R(m^{\star},\varphi)$ in terms of $\bE U(m^{\star},r,\varphi)$:
\begin{eqnarray*}
\lefteqn{\bE U(m^{\star},r,\varphi)} \\
   & = & \bE \left[ - \sum_{k=1}^{m^{\star}} \hat\lambda_k^{-2} \hat r_k^2 + \frac{\log^2 n}{n} \sum_{k=1}^{m^{\star}} \hat\lambda_k^{-2} \hat\sigma_k^2 \right],\\
   & = & \bE \left[ - \sum_{k=1}^{m^{\star}} \lambda_k^{-2} r_k^2 + \frac{\log^2 n}{n} \sum_{k=1}^{m^{\star}} \hat\lambda_k^{-2} \sigma_k^2 \right] - \bE \left[ \sum_{k=1}^{m^{\star}} \lbrace \hat\lambda_k^{-2} \hat r_k^2 - \lambda_k^{-2} r_k^2 \rbrace \right] \\
   &  & \hspace{2cm} - \frac{\log^2 n}{n} \bE \left[ \sum_{k=1}^{m^{\star}} \hat\lambda_k^{-2} (\sigma_k^2 - \hat\sigma_k^2) \right], \\
   & = & \bE \left[ \sum_{k>m^{\star}} \varphi_k^2 + \frac{\log^2 n}{n} \sum_{k=1}^{m^{\star}} \hat\lambda_k^{-2} \sigma_k^2 \right] - \|\varphi \|^2 -\bE \left[ \sum_{k=1}^{m^{\star}} \lbrace \hat\lambda_k^{-2} \hat r_k^2 - \lambda_k^{-2} r_k^2 \rbrace \right] \\
   &  & \hspace{2cm} - \frac{\log^2 n}{n} \bE \left[ \sum_{k=1}^{m^{\star}} \hat\lambda_k^{-2} (\sigma_k^2 - \hat\sigma_k^2) \right]. 
\end{eqnarray*}
Hence,
\begin{eqnarray}
\bE \bar R(m^{\star},\varphi) & = & \bE U(m^{\star},r,\varphi) + \|\varphi \|^2 + \bE \left[ \sum_{k=1}^{m^{\star}} \lbrace \hat\lambda_k^{-2} \hat r_k^2 - \lambda_k^{-2} r_k^2 \rbrace \right] \nonumber \\
   & & \hspace{2cm} +\frac{\log^2 n}{n} \bE \left[ \sum_{k=1}^{m^{\star}} \hat\lambda_k^{-2} (\sigma_k^2 - \hat\sigma_k^2) \right].
\label{num}
\end{eqnarray}
Remark that: 
\begin{eqnarray*}
\lefteqn{\bE \left[ \sum_{k=1}^{m^{\star}} \lbrace \hat\lambda_k^{-2} \hat r_k^2 - \lambda_k^{-2} r_k^2 \rbrace \right]}\\
& = & \bE \left[ \sum_{k=1}^{m^{\star}} \hat\lambda_k^{-2}( \hat r_k^2 -  r_k^2) \right]+ \bE \left[ \sum_{k=1}^{m^{\star}} (\hat\lambda_k^{-2} - \lambda_k^{-2}) r_k^2\right],\\
& = & \bE \left[ \sum_{k=1}^{m^{\star}} \hat\lambda_k^{-2} \lbrace (\hat r_k - r_k)^2 +2(\hat r_k -r_k)r_k \rbrace \right] + \bE \left[ \sum_{k=1}^{m^{\star}} (\hat\lambda_k^{-2} - \lambda_k^{-2}) r_k^2\right].
\end{eqnarray*}
Using simple algebra:
\begin{eqnarray*}
\lefteqn{\bE \sum_{k=1}^{m^{\star}} \hat\lambda_k^{-2} (\hat r_k - r_k) r_k} \\
& = & \bE \sum_{k=1}^{m^{\star}} \lambda_k^{-2} (\hat r_k - r_k) r_k + \bE \sum_{k=1}^{m^{\star}} (\hat\lambda_k^{-2}- \lambda_k^{-2}) (\hat r_k - r_k) r_k ,\\
& = & \bE \sum_{k=1}^{m^{\star}} \lambda_k^{-2} (\hat r_k - r_k) r_k + \bE \sum_{k=1}^{m^{\star}} (\hat\lambda_k^{-1}- \lambda_k^{-1}) r_k  (\hat\lambda_k^{-1}+ \lambda_k^{-1})(\hat r_k - r_k),\\
& \leq & \bE \sum_{k=1}^{m^{\star}} \lambda_k^{-2} (\hat r_k - r_k) r_k + \bE \sum_{k=1}^{m^{\star}} (\hat\lambda_k^{-1}- \lambda_k^{-1})^2 r_k^2 + C\bE \sum_{k=1}^{m^{\star}} \hat\lambda_k^{-2} (\hat r_k -r_k)^2. 
\end{eqnarray*}
Hence,
\begin{eqnarray*}
\bE \left[ \sum_{k=1}^{m^{\star}} \lbrace \hat\lambda_k^{-2} \hat r_k^2 - \lambda_k^{-2} r_k^2 \rbrace \right]
& \leq & C \bE \left[ \sum_{k=1}^{m^{\star}} \hat\lambda_k^{-2}  (\hat r_k - r_k)^2 \right] + 2 \bE \left[ \sum_{k=1}^{m^{\star}} \lambda_k^{-2} (\hat r_k -r_k)r_k \right] \\
& &  +\bE \left[ \sum_{k=1}^{m^{\star}} (\hat\lambda_k^{-2} - \lambda_k^{-2}) r_k^2\right] + \bE \sum_{k=1}^{m^{\star}} \left( \frac{\lambda_k}{\hat\lambda_k}-1 \right)^2 \varphi_k^2.
\end{eqnarray*}
Using Lemmata 3.4, \ref{risquevpbruitee} and (\ref{numT2}), we obtain, for all $1>\gamma>0$ and $K>1$:
\begin{eqnarray}
\lefteqn{\bE \left[ \sum_{k=1}^{m^{\star}} \lbrace \hat\lambda_k^{-2} \hat r_k^2 - \lambda_k^{-2} r_k^2 \rbrace \right]} \nonumber \\
& \leq & \left( 2 \gamma^{-1} \log^K n + C\gamma^{-1}\log^{3/2} n +\gamma \right). \frac{1}{n} \bE \left[ \sum_{k=1}^{m^{\star}} \hat\lambda_k^{-2} \sigma_k^2 \right] \label{blamstar} \\
& &   + \gamma^{-1} R(m_0,\varphi) + \gamma \bE \left[ \sum_{k>m^{\star}} \varphi_k^2 \right] + \Omega +  C\gamma^{-1} N^{2t +1 }e^{-\log^K n}+ \frac{C}{n}\left( \frac{\log^2(n).\|\varphi\|^2}{\gamma} \right)^{2t}. \nonumber
\end{eqnarray}
Remark that this result can be obtained for all $\bar m$ measurable with respect to the sample $(X_i,Y_i,W_i)_{i=1\dots n}$. Then, from (\ref{num}) and Lemma 3.6,
\begin{eqnarray*}
\lefteqn{\bE \bar R(m^{\star},\varphi)}\\
& \leq & \bE U(m^{\star},r,\varphi) + \|\varphi \|^2+ \left( 2\gamma^{-1}\log^K n + C\gamma^{-1} \log^{3/2} n + C\frac{\log^2 n}{n^{1/2}}\right) \frac{1}{n} \bE \left[\sum_{k=1}^{m^{\star}} \hat\lambda_k^{-2} \sigma_k^2 \right]\\
& & \hspace{-0.5cm}  + \gamma^{-1} R(m_0,\varphi) + \gamma \bE \left[ \sum_{k> m^{\star}} \varphi_k^2 \right] + C\gamma^{-1} N^{2t+1} e^{-\log^K n}+\Omega + \frac{C}{n}\left( \frac{\log^2(n).\|\varphi\|^2}{\gamma} \right)^{2t}.
\end{eqnarray*}
which can be rewritten:
\begin{eqnarray}
(1-\rho(\gamma,K,n))\bE \bar R(m^{\star},\varphi)
& \leq & \bE U(m^{\star},r,\varphi) + \|\varphi \|^2 \nonumber \\
& & \hspace{-4cm} + 2\gamma^{-1} R(m_0,\varphi) + C\gamma^{-1} N^{2t+1} e^{-\log^K n}  + \Omega + \frac{C}{n} \left( \frac{\log^2(n).\|\varphi\|^2}{\gamma} \right)^{2t},
\label{Etape2}
\end{eqnarray}
with 
$$ \rho(\gamma,K,n) = 2\gamma^{-1}\log^{K-2} n+ \frac{C}{n^{1/2}}+\log^{-1/2}n +\gamma.$$
The third step of our proof can be easily derived from the definition of $m^{\star}$ and leads to the following result: 
\begin{eqnarray}
(1-\rho(\gamma,K,n))\bE \bar R(m^{\star},\varphi) & \leq & \bE U(m_1,r,\varphi) + \|\varphi \|^2 + 2\gamma^{-1} R(m_0,\varphi) \nonumber \\
& &   +C\gamma^{-1} N^{2t+1} e^{-\log^K n}  + \Omega + \frac{C}{n}\left( \frac{\log^2(n).\|\varphi\|^2}{\gamma} \right)^{2t},
\label{Etape3}
\end{eqnarray}
where $m_1$ is defined in \eqref{eq:m1} and denotes the oracle in the family $\lbrace 1,\dots, M \rbrace$. In order to conclude the proof, we have to compute $\bE U(m_1,r,\varphi) + \|\varphi \|^2$. In a first time, remark that:
\begin{eqnarray*}
\bE U(m_1,r,\varphi) + \|\varphi \|^2
& = & \bE\left[ - \sum_{k=1}^{m_1} \hat\lambda_k^{-2} \hat r_k^2 + \frac{\log^2 n}{n} \sum_{k=1}^{m_1} \hat\lambda_k^{-2} \hat\sigma_k^2 \right] + \|\varphi \|^2,\\
& = & \bE \left[ - \sum_{k=1}^{m_1} \lambda_k^2 r_k^2 \right] + \|\varphi\|^2 + \frac{\log^2 n}{n} \bE \left[ \sum_{k=1}^{m_1} \hat \lambda_k^{-2} \sigma_k^2 \right] \\
& & + \bE\left[ \sum_{k=1}^{m_1} (\lambda_k^{-2} r_k^2 - \hat\lambda_k^{-2} \hat r_k^2)\right] + \frac{\log^2 n}{n} \bE \left[ \sum_{k=1}^{m_1} (\hat\lambda_k^{-2} \hat\sigma_k^2 - \hat\lambda_k^{-2} \sigma_k^2 ) \right].
\end{eqnarray*}
Hence,
\begin{eqnarray*}
\lefteqn{\bE U(m_1,r,\varphi) + \|\varphi \|^2}\\
& = & \bE \left[ \sum_{k>{m_1}} \varphi_k^2 + \frac{\log^2 n}{n} \sum_{k=1}^{m_1} \hat \lambda_k^{-2} \sigma_k^2 \right]+ \bE\left[ \sum_{k=1}^{m_1} (\lambda_k^{-2} r_k^2 - \hat\lambda_k^{-2} \hat r_k^2)\right] \\
& &  \hspace{1cm} + \frac{\log^2 n}{n} \bE \left[ \sum_{k=1}^{m_1} (\hat\lambda_k^{-2} \hat\sigma_k^2 - \hat\lambda_k^{-2} \sigma_k^2 ) \right],\\
& = & \bE \bar R(m_1,\varphi) + F_1 + F_2.
\end{eqnarray*}
The same bound as (\ref{blamstar}) occurs for $F_1$. By the same way, using Lemma 3.6:
\begin{eqnarray*}
F_2 & = & \frac{\log^2 n}{n} \bE \left[ \sum_{k=1}^{m_1} (\hat\lambda_k^{-2} \hat\sigma_k^2 - \lambda_k^{-2} \sigma_k^2 ) \right],\\
    & \leq & C\frac{\log n}{n^{3/2}}. \bE \left[ \sum_{k=1}^{m_1} \hat\lambda_k^{-2} \sigma_k^2 \right] + \frac{1}{n} \bE \sum_{k=1}^{m_1} \hat\lambda_k^{-2}(r_k^2 - \hat r_k^2) + C e^{-\log^2 n}.
\end{eqnarray*}
Therefore, for all $K\geq 1$,
\begin{eqnarray}
\bE U(m_1,r,\varphi) + \|\varphi \|^2
& \leq & \left( 1+C \log^{K-2} n + \frac{C\log^{-1} n}{\sqrt{n}} \right)\bE \bar R(m_1,\varphi) + R(m_0,\varphi) \nonumber\\
& &  + C \gamma^{-1} N^{2t +1} e^{-\log^K n} + \frac{C}{n} \left( \frac{\log^2(n).\|\varphi\|^2}{\gamma} \right)^{2t} + \Omega.
\label{lastnum}
\end{eqnarray}
Using (\ref{Etape3}) and (\ref{lastnum}), we eventually obtain:
\begin{eqnarray*}
\lefteqn{(1-\rho(\gamma,K,N)) \bE \bar R(m^{\star},\varphi)}\\
& \leq & \left( 1+ \log^{K-2} n + \frac{C\log^{-1} n}{\sqrt{n}} \right) \bE \bar R(m_1,\varphi)+C\gamma^{-1} \bE R(m_0,\varphi) \\
 & & + C \gamma^{-1} N^{2t +1}  e^{-\log^K n} + \frac{C}{n} \left( \frac{\log^2(n).\|\varphi\|^2}{\gamma} \right)^{2t} + \Omega,\\
& \leq & C \log^2 (n). \bE R(m_1,\varphi) + C\gamma^{-1} \bE R(m_0,\varphi) + C\gamma^{-1} N^{2t+1} e^{-\log^K n} + \frac{C}{n} \left( \frac{\log^2(n).\|\varphi\|^2}{\gamma} \right)^{2t} +\Omega,\\
& \leq & C \log^2 (n). R(m_0,\varphi) + \log^2(n). \Gamma(\varphi) + \frac{C}{n} \left( \frac{\log^2(n).\|\varphi\|^2}{\gamma} \right)^{2t} + \Omega,
\end{eqnarray*}
for some positive constant $C$, where $\Gamma(\varphi)$ is introduced in Theorem 1. With an appropriate choice of $K$, this leads to:
\begin{eqnarray*}
\lefteqn{\bE\| \varphi^{\star} - \varphi \|^2}\\
& \leq  & C \log^2 (n). R(m_1,\varphi) + \frac{C}{n} \left( \frac{\log^2 (n). \|\varphi\|^2}{\gamma} \right)^{2t} + \Omega + \log^2 (n). \Gamma(\varphi).
\end{eqnarray*}
\begin{flushright}
$\Box$
\end{flushright}

\noindent
\textsc{Proof of Corollay~\ref{larate}} We start by recalling the oracle inequality obtained for the estimator $\vphi^\star$.
\begin{eqnarray*}
\mathbb{E} \| \varphi^{\star} - \varphi \|^2 & \leq  & C_0 \log^2 (n). \left[ \inf_m R(m,\varphi) \right]+ \frac{C_1}{n} \left( \log (n). \|\varphi \|^2 \right)^{2\beta} \\
& & \hspace{3cm} +\Omega + \log^2 (n) .\Gamma(\varphi),
\end{eqnarray*}
We have to bound the risk under the regularity condition and the extra term $\log^2(n) \Gamma(\vphi)$.
Recall that the risk is given by
$$R(m,\varphi) = \sum_{k>m} \varphi_k^2 + \frac{\log^2 n}{n} \sum_{k=1}^m \lambda_k^{-2} \sigma_k^2.$$
Hence under \eqref{eq:Sobol}, we obtain both upper bounds for two constants $C_1$ and $C_2$
$$ \sum_{k>m} \varphi_k^2 \leq m^{-2s} C_1,$$
$$\frac{\log^2 n}{n} \sum_{k=1}^m \lambda_k^{-2} \sigma_k^2 \leq C_2 \frac{\log^2 n}{n} \sigma_U^2 m^{2t+1}.$$
An optimal choice is given by $m=[( n / \log n)^{\frac{1}{1+2s+2t}}]$, leading to the desired rate of convergence.\\
\indent Now consider the remainder term $\Gamma(\varphi)$. Under Assumption [IP], $M_0 \geq [n^{1/2s}/ \log^2 n]$, but since $m_0=[n^{\frac{1}{1+2s+2t}}]$ we get clearly that  $m_0 \leq M_0$, which entails that $\Gamma(\vphi)=0$.

\bibliographystyle{alpha}
\bibliography{Ivbiblio}
\noindent Jean-Michel {\sc Loubes}\hfill  Cl\'ement {\sc Marteau}\\
\noindent Equipe de probabilites et statistique,\hfill Insistut de 
Math\'{e}matique de Toulouse,\\ Insistut de 
Math\'{e}matique de Toulouse,\hfill  INSA D\'{e}partement GMM.\\
UMR5219, Universit\'e de Toulouse,\hfill \\
\noindent 31000 Toulouse  {\sc France}\hfill 
$\,$
\\
\noindent Jean-Michel.Loubes@math.univ-toulouse.fr \hfill  Clement.Marteau@insa-toulouse.fr\\
\end{document}

%% file: macros.tex
\newcommand{\DOIFPDF}[2]{\ifx\pdfoutput\undefined #2\else#1\fi}

\newcommand{\EM}{\ensuremath}

\newcommand{\PDFABLE}[2]{%
\newif\ifpdf%
\ifx\pdfoutput\undefined\pdffalse%
\else\pdftrue\pdfoutput=1\pdfcompresslevel=9\fi%
\ifpdf%
 \usepackage#1
 \usepackage[pdftex,%
             a4paper,%
             colorlinks,%
             citecolor=blue,%
             pagebackref,%
             plainpages=false]{hyperref}%
\else%
 \usepackage#2
 \usepackage{url}
\fi}

\makeatletter
\newcommand{\@THMSTYLES}{%
  \newtheoremstyle{bodyrm}
  {3pt}
  {3pt}
  {}
  {}
  {\bfseries\sffamily}
  {.}
  { }
  {}
  \newtheoremstyle{bodyit}
  {3pt}
  {3pt}
  {\itshape}
  {}
  {\bfseries\sffamily}
  {.}
  { }
  {}
}
\newcommand{\THMEN}{%
  \@THMSTYLES
  \theoremstyle{bodyit}
  \newtheorem{thm}{Theorem}[section]%
  \newtheorem{cor}[thm]{Corollary}%
  \newtheorem{prop}[thm]{Proposition}%
  \newtheorem{lem}[thm]{Lemma}%
  \theoremstyle{bodyrm}%
  \newtheorem{defi}[thm]{Definition}%
  \newtheorem{xpl}[thm]{Example}%
  \newtheorem{exo}[thm]{Exercise}%
  \newtheorem{hyp}[thm]{Hypothesis}%
  \newtheorem{eur}[thm]{Heuristics}%
  \newtheorem{pro}[thm]{Problem}%
  \newtheorem{rem}[thm]{Remark}%
  \newtheorem{prp}[thm]{Property}%
}
\newcommand{\THMFR}{%
  \@THMSTYLES
  \theoremstyle{bodyit}
  \newtheorem{thm}{Théorème}[section]%
  \newtheorem{cor}[thm]{Corollaire}%
  \newtheorem{lem}[thm]{Lemma}%
  \theoremstyle{bodyrm}%
  \newtheorem{rem}[thm]{Remarque}%
}
\makeatother

\makeatletter
\newcommand{\SMALLSECS}{%
 \renewcommand{\section}{\@startsection%
  {section}
  {1}
  {0em}
  {\baselineskip}
  {0.5\baselineskip}
  {\normalfont\large\bfseries}}
 \renewcommand{\subsection}{\@startsection%
  {subsection}
  {2}
  {0em}
  {\baselineskip}
  {0.25\baselineskip}
  {\normalfont\bfseries}}
}
\makeatother

\makeatletter
\providecommand{\timenow}{\@tempcnta\time
\@tempcntb\@tempcnta
\divide\@tempcntb60
\ifnum10>\@tempcntb0\fi\number\@tempcntb
\multiply\@tempcntb60
\advance\@tempcnta-\@tempcntb
:\ifnum10>\@tempcnta0\fi\number\@tempcnta}
\makeatother

\makeatletter
\newcommand{\versiondetravail}{%
 \renewcommand{\@evenfoot}{%
 \hfil{\tiny\texttt{%
   Version préliminaire, compilée le \today{} à \timenow.}\hfill}}%
 \renewcommand{\@oddfoot}{\@evenfoot}%
}
\makeatother

\newcommand{\dN}{\EM{\mathbb{N}}}

\newcommand{\bE}{\EM{\mathbf{E}}}

\newcommand{\bR}{\EM{\mathbf{R}}}

\newcommand{\la}{\lambda}

\newcommand{\vphi}{\varphi}

\newcommand{\Det}[1]{\mathrm{Det}\,}

\renewcommand{\leq}{\leqslant}
\renewcommand{\geq}{\geqslant}

